\documentclass[12pt, final]{l4dc2023}

\usepackage{times}
\newcommand{\norm}[1]{\left\lVert#1\right\rVert}

\newtheorem{assumption}[theorem]{Assumption}

\title[SLS for Nonlinear Systems]{Designing System Level Synthesis Controllers for Nonlinear Systems with Stability Guarantees}

\author{\Name{Lauren Conger} \Email{lconger@caltech.edu}\\
 \Name{Sydney Vernon} \Email{svernon@caltech.edu}\\
 \Name{Eric Mazumdar} \Email{mazumdar@caltech.edu}\\
 \addr 1200 E. California Blvd., Pasadena, CA 91125}

\begin{document}

\maketitle

\begin{abstract}%
We introduce a method for controlling systems with nonlinear dynamics and full actuation by approximating the dynamics with polynomials and applying a system level synthesis controller. We show how to optimize over this class of controllers using a neural network while maintaining stability guarantees, without requiring a Lyapunov function. We give bounds for the domain over which the use of the class of controllers preserves stability and gives bounds on the control costs incurred by optimized controllers. We then numerically validate our approach and show improved performance compared with feedback linearization--- suggesting that the SLS controllers are able to take advantage of nonlinearities in the dynamics while guaranteeing stability.
\end{abstract}


\section{Introduction}
Many physical systems systems are modeled by nonlinear dynamics, and these systems need to be controlled. One existing option is approximating nonlinear systems systems with polynomials and using sum of squares to construct controllers \citep{strasser_data-driven_2021}. Other successful techniques include model predictive control (MPC), learning dynamics or controllers with neural networks \citep{shi_neural_2019,shi_neural-swarm2_2022,furieri_neural_2022}, and control Lyapunov functions (CLFs)\citep{ames_rapidly_2014}. Tools for nonlinear control are less developed than those for linear control because of the breadth and complexity of various nonlinear systems, so locally linearizing the system along a trajectory or lifting to linear dynamics \citep{lovchakov_multidimensional_2020} allows engineers to utilize linear system theory.  However, sometimes the accuracy of the model must be better than is feasible with linearized models, as in the case of modeling braking of wheeled vehicles \citep{vataeva_parametric_2019}, motivating us to pursue a controller in the nonlinear setting.

We choose a polynomial form of the dynamics because many nonlinear systems can be lifted into this form \citep{qian_lift_2020,mauroy_linear_2016} and polynomial dynamics and state feedback controller can be used to ensure stable periodic orbits \citep{yuno_realization_2014}. Additionally, polynomials have been successfully used as a basis for identifying nonlinear Hammerstein systems \citep{hammar_fractional_2015}.
We offer a method for discrete-time control for systems with nonlinear dynamics, where the system is fully actuated, with a polynomial approximation of the dynamics and a nonlinear system level synthesis (SLS) controller resulting in input-to-state stable dynamics with respect to external disturbances. The majority of work uses SLS, the basis of which is described in \citep{anderson_system_2019}, and \citep{ho_system_2020} describes conditions that transfer functions must satisfy in the nonlinear SLS setting.

\subsection{Contributions}
We propose a new method for controlling nonlinear systems in which we approximate nonlinear dynamics with polynomial terms, and design a system level synthesis (SLS) controller for the polynomial approximation. We show that using the controller on the original system results in stable closed-loop dynamics with respect to disturbances. Previous work in SLS  focused on controlling polynomial systems \citep{conger_nonlinear_2022}, and we show how to extend this technique to general nonlinear systems, and how such a controller can be represented with function approximators like neural networks---while preserving stability--- and then optimized to minimize a control objective. Furthermore, we give bounds on the control cost incurred by the optimized controllers in terms of the degree of the polynomial approximation and the effective domain of approximation. We empirically observe that the optimized nonlinear controllers have improved performance when compared with feedback linearization in a numerical example. 

\paragraph{Organization:} The remainder of the paper is organized as follows. We give a brief overview of related work and in Section~\ref{sec:prelims} we give details about the class of nonlinear systems that we consider and background for the polynomial SLS controller. In Section~\ref{sec:approx}, we derive a class of SLS controllers that are based on a polynomial approximation of the nonlinear dynamics that still have stability guarantees for the true nonlinear system. In Section~\ref{sec:optimize} we discuss how to parameterize the controller by a neural network and optimize the controller. We give a bound on the control cost incurred by a properly optimized controller as a function of the disturbance size and approximation order. Lastly, in Section~\ref{sec:example},  we illustrate our method on a point mass system that requires the neural network to parameterize and optimize 14 controller functions, comparing our results with an unoptimized controller and feedback linearization. We conclude with a brief discussion of future work.

\subsection{Related Work}
Previous work \citep{conger_nonlinear_2022} proved stability of the nonlinear SLS controller \eqref{eq:SLS_controller} for polynomial systems, and we extend this work by giving a performance bound on the control cost for $k$-times differentiable functions, a much broader class of nonlinear dynamics than just polynomial functions. The work by \citep{furieri_neural_2022} offers a solution for learning stabilizing controllers, similar to our work, but requires a stabilizing controller initially; our work provides that stabilizing controller. The criteria for general operator stability and robustness to model errors is developed in \cite{ho_system_2020} which inspired our ideas for contributing specific operators that satisfy that criteria. 

\section{Preliminaries}\label{sec:prelims}
We begin by describing the class of nonlinear systems which we can estimate with polynomial dynamics. We focus on nonlinear dynamics of the form:
\begin{align}\label{eq:nonlinear_dynamics}
    x_{t+1} = f(x_t) + u_t + w_t,
\end{align}
where $w_t\in\mathbb{R}^n$ is a disturbance, $x_t\in\mathbb{R}^n$ is the state, and $u_t\in\mathbb{R}^n$ is the input. In this work we focus on the case where the system has an equilibrium point at $x^*$ and the system is being perturbed by noise near $x^*$. We remark that the equilibrium does not need to be stable. Given this setup, we make the following assumption on the dynamics.
\begin{assumption}
Assume that the nonlinear function $f\in C^{\infty}$ has bounded derivatives up to order $k$: \[ \partial^{p} f(x_t) \leq M,\] for all $|p|\leq k$. Further assume that the noise is bounded so that $\norm{w_t}_\infty \leq W$.
\end{assumption}
The assumption on bounded noise allows us to treat the case with adversarial disturbances, while the assumption on bounded partial derivatives of up to order $k$ is necessary to give our stability and boundedness guarantees for the SLS controllers. This is a generalization of the common Lipchitzness assumptions in nonlinear control, where in essence we show that one can leverage higher order smoothness to design stabilizing controllers. 

\subsection{Polynomial SLS Controller}
Key to our controller design is a polynomial approximation of the nonlinear dynamics ~\eqref{eq:nonlinear_dynamics} and previous work on designing SLS controllers for polynomial systems~\citep{conger_nonlinear_2022}. Indeed, for a polynomial system of the form:
\begin{equation}\label{eq:poly_vector_dynamics}
    \begin{split}
        x_{t+1} = \sum_{j=1}^n A_j x_t^{\otimes j} + u_t + w_t,
    \end{split}
\end{equation}
where $A_j \in \mathbb{R}^{n^j \times n^j}$ and $\otimes$ denotes the column-wise Kronecker product of $x_t$ with itself $j-1$ times, a stabilizing finite impulse response (FIR) controller with time horizon $T$ is given by:

\begin{align}\label{eq:SLS_controller}
    u_t =& - \sum_{m=0}^{T-1} \sum_{j=1}^{c_m} \alpha_j^{(m)} G_j^{(m)}(w_{t:t-m}) -\sum_{j=1}^{c_T} G_j^{(T)}(w_{t:t-T}),
\end{align}
where $w_{t:t-m}$ denotes disturbances from time $t-m$ to time $t$, where $w\in\mathbb{R}^n$ \citep{conger_nonlinear_2022}. The functions $G_j^{(m)}$ are monomials of entries in the disturbance vectors, that is, $G_j^{(m)}(w_{t:t-m}) = d_j^{(m)}w_q(i)w_p(j)\hdots w_s(k)$ where each $w_q(i)$ denotes the $i$th entry of $q$th $w$ vector raised to some power that is less than or equal to $m$, and $d_j^{(m)}$ is a coefficient from the dynamics. The $c_m$ constant denotes the number of monomials that depend on the $w_m$ terms multiplied with any more recent disturbances, and we say $c = \max_{m\in(0,T-1)} c_m$, which is the maximum number of terms in the inner sum. Lastly, each $\alpha_j^{(m)}\in[0,1]$ weights each monomial term.

This controller results in the state $x_t$ as a function of $w_{t-1:t-T}$
\begin{align}\label{eq:SLS_states}
    x_t =& \sum_{k=0}^{T-1} \sum_{j=1}^{c_m}(1-\alpha_j^{(m)})G_j^{(m)}(w_{t-1:t-1-k}) + w_{t-1}.
\end{align}
which in turn implies that the state is uniformly bounded across time. This is a crucial step in our later analysis.

\section{Polynomial Approximation of Nonlinear Dynamics}~\label{sec:approx}
Given the setup described above, we now show that nonlinear systems of the form \eqref{eq:nonlinear_dynamics} can be stabilized by SLS controllers designed based on polynomial approximations to the dynamics. 

To begin, we assume that we have access to a polynomial approximation to the nonlinear dynamics around the (potentially unstable) equilibrium. We assume that the approximation is derived using a Taylor expansion, rather than a Chebyshev approximation, because we are stabilizing to a single point and would like the highest accuracy of the estimation at that point. This allows us to rewrite the nonlinear dynamics as:
\begin{align}\label{eq:closed_loop_dyn_with_approx}
    x_{t+1} = \sum_{j=1}^k H_j x_t^{\otimes j} + u_t + e(x_t) + w_t
\end{align}
where $e(x_t)$ is the error from the Taylor series approximation, and $e(0)=0$.
Given this form,  we proceed by designing a SLS controller tailored to this polynomial approximation and show that it still has stability guarantees for the true system. Indeed, in the next section we illustrate the order-of-magnitude tradeoffs between the order of the polynomial approximation, domain over which a bounded control is required for control, and state dimension of the system. Intuitively, we expect a higher order of approximation to result in a larger region over which we can stabilize with bounded control, and that a larger state dimension will shrink this region, which we show via a specific relationship.

\subsection{Stability Guarantee Under Polynomial Approximation}
We now analyze the stability of the SLS controller \eqref{eq:SLS_controller}  based off of the polynomial approximation to the nonlinear system. In particular, we prove that the controller is stabilizing over a certain domain around the equilibrium. To do so, we use some definitions that facilitate the proof. We denote $\alpha = \min_{j,m} \alpha_j^{(m)}$, and note that $G_j^{(m)}(w_{t:t-m}) < l \norm{ w_{t:t-m}}_\infty$ for some $l>0$.  Note that the linear bound on $G_j^{(m)}$ is possible because $w_t$ is bounded for all $t$. Lastly, define input-to-state stability (ISS) with respect to the disturbance as a guarantee that $\norm{x_t} \leq \gamma(\norm{W})$ for the initial condition $\norm{x_0}=0$, where $\gamma$ is some class $\mathcal{K}$ function.

\begin{theorem}\label{thm:bounded_state}
The system is ISS using approximate dynamics with the SLS controller \eqref{eq:SLS_controller} for $k$ large enough and all $\alpha>0 $ such that 
\begin{align*}
    \frac{MW^k}{(k+1)!}+lc(1-\alpha) < 1.
\end{align*}
\end{theorem}
This gives a continuum of $\alpha$ over which to optimize while guaranteeing stability. Considering the case where $\alpha_j^{(m)}=1$, which means disturbances are canceled as quickly as possible, we are left with the condition that $\frac{MW^k}{(k+1)!}\leq 1$
which illustrates how $k$ must be sufficiently large. Note that this bound is not tight, and larger values of $\alpha$ than the ones satisfying the lemma result in stable behavior in simulation. 

\begin{proof}
We use the Lagrange form of the remainder for the Taylor expansion to bound the error term $e(W) \leq \frac{MW^{k+1}}{(k+1)!}$. To show stability we will show that the state remains bounded given bounded noise and $x_0=0$. We use a recursion relation to show how the state evolves and bound the final state. Starting with the first state, we have that $\norm{x_1} \leq W$ since the first state was zero and was hit with a disturbance of size at most $W$. The second state, computed from \eqref{eq:closed_loop_dyn_with_approx}, will also have the error from the approximate dynamics and the partially cancelled terms left by the controller
\begin{align*}
    \norm{x_{2}} \leq (1-\alpha)lc W +  \left(\frac{MW^k}{(k+1)!}+1 \right) W.
\end{align*}
Computing this again for the next state, $x_3$, we have
\begin{align*}
    \norm{x_3} \leq W\bigg[ (1-\alpha)^2(lc)^2 + (1-\alpha)lc\left(2\frac{MW^k}{(k+1)!}+1\right) + \left(\left(\frac{MW^{k+1}}{(k+1)!}\right)^2+\frac{MW^k}{(k+1)!}+1\right) \bigg]
\end{align*}
and we can continue computing an upper bound on $\norm{x_t}$ as 
\begin{align*}
    \norm{x_t} \leq W\sum_{\tau=0}^{t-1} \sum_{i=0}^{\tau} \binom{\tau}{i} \left(\frac{MW^k}{(k+1)!}\right)^i  (1-\alpha)^{\tau-i}(lc)^{\tau-i}.
\end{align*}
Simplifying the above summation yields $\norm{x_t} \leq W\frac{b^t-1}{b-1}$ where $b=\frac{MW^k}{(k+1)!}+lc(1-\alpha)$. From here we deduce that in order for the sum to be finite, the system must satisfy
\begin{align*}
    \frac{MW^k}{(k+1)!}+lc(1-\alpha) < 1.
\end{align*}
\end{proof}

\section{Optimizing the Controller}~\label{sec:optimize}
Given the set of stabilizing controllers in the previous section, we now give a method for parametrizing the controllers and optimizing their performance on a control objective. In particular we allow the coefficients $\alpha_j^{(k)}$ to be functions of the disturbances $w$ and use a function approximator (e.g., a neural network) to learn the coefficients $\alpha_j^{(k)}$ as a function of the past disturbances. 
Crucially, we show that an optimized controller results in a bounded control effort. 

We remark that previous work treated  $\alpha_j^{(k)}$  as constant coefficients. However, allowing them to be functions broadens the class of controllers that \eqref{eq:SLS_controller} parameterizes, leading to better performance. Let $H_\theta:\mathbb{R}^{nk} \rightarrow (\alpha,1)^k$ represent the function parameterized by a neural network with weights $\theta$. We assume without loss of generality that $H_\theta$ is a class of functions such that there exists a $\theta \in \Theta$ such that $H_\theta \equiv 1$. We minimize the quadratic control cost 
\begin{align}\label{eq:optimization_problem}
    \min_{\theta} &\sum_{t=1}^{N_T}  u_t^\top R u_t \\ \nonumber
    \alpha_j^{(m)} &= H_{\theta}(w_{t-1:t-1-m}) \\ \nonumber
    u_{t+1} &= -\sum_{m=0}^{T-1} \sum_{j=1}^{c_m} \alpha_j^{(m)} G_j^{(m)}(w_{t:t-m}) -\sum_{j=1}^{c_T} G_j^{(T)}(w_{t:t-T}) \\ \nonumber
     x_{t+1} &= f(x_t) + u_t + w_t \\ \nonumber
     x_0 &= 0
\end{align}
for $R \succ 0$ and $N_T$ the number of time steps for which we run the training disturbance series. Note that because the controller forces the system to have a finite impulse response, a cost on the state is not necessary for stabilization. The NN trains on random noise sequences that are longer than the FIR horizon and is penalized by the quadratic cost function, with the states and control inputs computed in the loss function. Note that the true dynamics $f$ should be used to generate the next state, and the approximated dynamics are used to compute the controller. The controller is illustrated in \ref{fig:nn_diagrams} for a FIR horizon of $T=3$.

\subsection{Tradeoffs}
Once the controller has been optimized, the $\alpha_j^{(m)}$ coefficients will be tuned so disturbances that decay are propagating through the dynamics and canceled at the end of the FIR horizon, and disturbances that are not decaying are cancelled immediately. This will have a lower cost than a controller that cancels all of the disturbances immediately, because cancelling a decayed disturbance is less costly than cancelling the larger disturbance. 
We upper bound the cost of the optimized controller by the control cost of setting all $\alpha_j^{(m)}=1$. We show that cost is bounded by a constant, and illustrate tradeoffs for the size of the domain over which the approximation is valid, the order (power) of the polynomial approximation, and the input size required to stabilize the system. We denote the desired domain over which stability is guaranteed is $h$ so that $\norm{x_t}_1\leq h$ which a design decision that can be satisfied by choosing the number of polynomial terms $k$ in accordance with Theorem~\ref{thm:bounded_state}.
\begin{theorem}
The controller $u_t^1$ with all $\alpha_j^{(m)}=1$ subject to disturbances $\{w_i\}_{i=0}^T$ results in:
\[\|u_t^1\|_2 \le M \left( \frac{W d(W^k d^k - 1)}{W d-1} + \frac{W (1-W^k)}{W-1} \right)=U^1 \quad \forall t,\]  
where $d=\frac{M(nh)^k}{W(k+1)!}+1$.
\end{theorem}

Expanding this upper bound for $u_t^1$ shows that $\norm{u_t^1}_2 \leq M \left( \frac{(Wd)^{k+1}}{W d-1} -\frac{W^{k+1}}{W-1} -\frac{W d}{W d-1}+ \frac{W}{W-1} \right)$. Ignoring the last two terms, which are constant as $W$ becomes large, we notice that $\frac{(Wd)^{k+1}}{W d-1} -\frac{W^{k+1}}{W-1}$ has the same behavior as $(Wd)^k -W^k$. Examining how $(Wd)^k -W^k$ evolves for increasing $k$, we see that $(Wd)^k -W^k =\left( \frac{M(nh)^{k+1}}{(k+1)!}+W \right)^k - W^k$ goes to zero as $k$ goes to infinity, since $\frac{M(nh)^{k+1}}{(k+1)!}$ goes to zero. We see that a higher order approximation, $k$, makes the error term smaller for sufficiently large $k$, while increasing the domain over which we expect the approximation bound to hold, $h$, or the state dimension, $n$, makes the error term larger. The following corollary follows directly from the theorem because we include in the class of functions $H_\theta$ the possibility that $H_\theta \equiv 1$, which gives a cost greater than or equal to the optimal solution to \eqref{eq:optimization_problem}.

\begin{corollary}
The optimized controller with a polynomial approximation of order $k$, bounded disturbances, and desired stability radius $h$ therefore has a control cost at most
\[L(H_\theta,\{w_i\}_{i=0}^T) \le rTU^1\]  
where $r=\sigma_{max}(R)$ and $L$ is the cost of the optimization of \eqref{eq:optimization_problem}.
\end{corollary}
The corollary shows that the optimal control cost is finite, decreases for sufficiently large $k$,  and increases for increasing $h$.
\begin{proof}[Theorem 3]
We approximate $f$ using a Taylor expansion around the point to which we stabilize the system, $x^*$, and we assume without loss of generality $x^*=0$. The expansion takes the form
$p_k(x_t) = \sum_{j=1}^k H_j x_t^{\otimes j}$. Note that we do not include the $j=0$ term because it corresponds to $f(x^*)$ which equals zero.
The approximation error from the Taylor expansion is given by $e(x) = f(x) - p_k(x)$
where $p_k(x)$ is a polynomial approximation of order $k$. We can rewrite the dynamics as 
\begin{align}
    x_{t+1} = p_k(x_t) + u_t + w_t + e(x_t).
\end{align}
The error can be treated as a disturbance term. In the worst (performance) case, we do not allow the error to propagate and feedback linearize immediately. The total disturbance term from time $t$ will be determined at time step $t+1$, and can be counteracted at time step $t+1$, so $u_{t+1}=-p_k(x_{t+1})$ which can be expanded as
\begin{align*}
    u_{t+1} &= -p_k(w_t+e(x_t)) \\
        &= -\sum_{j=1}^k H_j (w_t + e(x_t))^{\otimes j} \\
        &= -\sum_{j=1}^k H_j \sum_{i=1}^j C_{j i} \left(w_t^{\otimes j-i} \otimes  e(x_t)^{\otimes i} \right)
\end{align*}
where the $C_{ij}$ matrices contain the binomial expansion coefficients. We take the 2-norm of $u_t$ to upper bound the control input, and then use the property that the 2-norm of the Kronecker product can be upper bounded by the product of the individual norms \citep{cochrane_nuclear_2020,kong_new_2018}. The norm of the control input is given by
\begin{align*}
    \norm{u_t}_2 &\leq \sum_{j=1}^k \norm{ H_j }_2 \sum_{i=1}^j \norm{ C_{j i} }_2 \norm{  w }_2^{j-i} \norm{ e(x_t) }_2^i. 
\end{align*}
We use an upper bound on the binomial matrices $\norm{C_{ij}}_2\leq j!\binom{j}{i}$, and the worse-case noise magnitude $W$ to further bound the expression for $\norm{u_t}_2$. We also denote $\norm{e(x_t)}_2 = e_x$ for ease of notation
\begin{align*}
    u_t &\leq \sum_{j=1}^k \norm{ H_j }_2 \sum_{i=1}^j \frac{j!^2}{(j-i)!i!} W^{j-i} e_x^i.
\end{align*}
Evaluating the inner sum in a symbolic solver gives the relation
\begin{align*}
    \sum_{i=1}^j \frac{j!}{(j-i)!i!}  \left(\frac{e_x}{W}\right)^i = \left(\frac{e_x}{W}+1\right)^j-1
\end{align*}
which we can use to simplify the bound $u_t$ to a single sum 
\begin{align*}
    u_t &\leq \sum_{j=1}^k \norm{ H_j }_2 j!\ W^j \left(\left(\frac{e_x}{W}+1\right)^j-1 \right).
\end{align*}
Now we use the definition of $\norm{H_j}_2$ as the matrix coefficients for the Taylor expansion of the function, which means each norm is upper-bounded by the bound of the $k+1$-th derivatives divided by a factorial, that is, $\norm{H_j}_2 \leq M/j!$. This leaves us with the sum
\begin{align*}
     u_t &\leq M\sum_{j=1}^k W^j \left(\left(\frac{e_x}{W}+1\right)^j-1 \right).
\end{align*}
Lastly, we can again evaluate this sum in a symbolic solver. Let $d=\frac{e_x}{W}+1$ for ease of notation, and then we have that $\norm{u_t}_2 \leq M \left( \frac{W d(W^k d^k - 1)}{W d-1} + \frac{W (1-W^k)}{W-1} \right)$. 
\end{proof}

\section{Example}\label{sec:example}
We illustrate the polynomial approximation on an example system. We approximate the system dynamics and compute the controllers. Then we parameterize the $\alpha_j^{(m)}(w_{t:t-m})$ using a NN and implement a quadratic cost function. After training the NN, we compare the total quadratic cost with feedback linearization under simulation scenarios in which the system is near its unstable equilibrium and continually battered by disturbances. We trained the network on Google Colab's freely available Intel(R) Xeon(R) CPU @ 2.20GHz using stochastic gradient descent.

\subsection{Neural Network Architecture}\label{subsec:NN_architecture}
\begin{figure}\label{fig:nn_diagrams}
    \centering
    \includegraphics[scale=0.5,trim=50 100 25 50, clip]{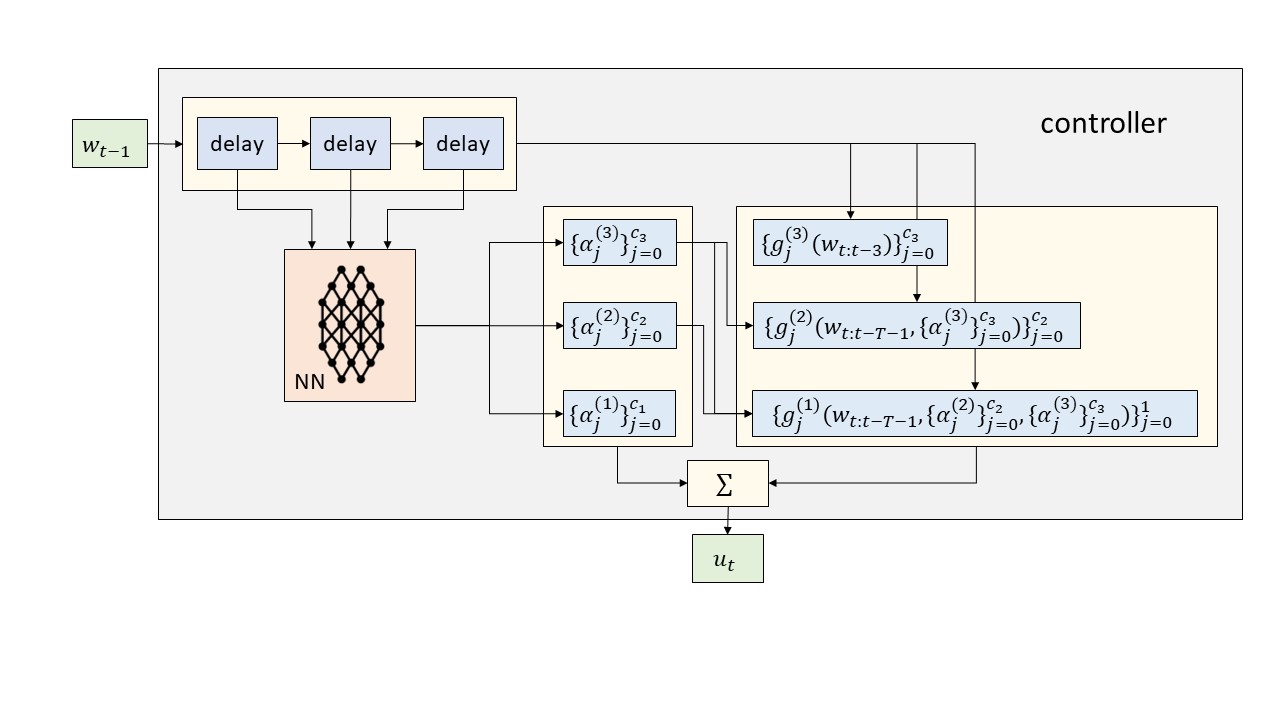}
    \caption{The controller diagram shows a FIR horizon of $T=3$. The NN generates values for each $\alpha$ coefficient at each time step, which are coefficients to the $G_j^{(m)}$ and used to compute the $G_j^{(m)}$ functions of more recent disturbances.}
    \label{fig:disturbances_to_alpha}
\end{figure}

The neural network accepts as input a list of the previous $T$ disturbances $w_{t:t-T}$. The disturbances pass through a NN with linear layers, ReLU activation functions, and dropout layers. The NN outputs the values assigned to each value of $\alpha_j^{(m)}$ corresponding to the most recent disturbance. Causality is enforced by computing each set of $\alpha_j^{(m)}$ at time $t$ using inputs $w_{t-1:t-1-m}$, so that disturbances are not from current or future time steps. The diagram in Figure~\ref{fig:nn_diagrams} illustrates how disturbances are used to compute coefficient weights.

\subsection{Point Mass}\label{subsec:point mass}
In this example we model a hovering point mass subject to disturbance. We seek to stabilize the point to a fixed vertical height, and minimize a linear quadratic cost function of the state and input. The discretized dynamics are given by
\begin{align}
    \ddot{z} = -g + \frac{1}{m}f_a(z,\dot{z}) + u + w
\end{align}
where $\ddot{z}$ denotes the acceleration, $g$ is the gravitational constant, $m$ is the mass of the point, and $f_a$ is a NN model for the acceleration dynamics. The model is taken from \citep[Eq 8]{shi_neural_2019}, which demonstrates a method for learning the dynamics of a system while flying and providing stability guarantees via Lipschitz constraints on the network layers.

\begin{figure}\label{fig:point_mass_cost}
    \centering
    \includegraphics[scale=0.5]{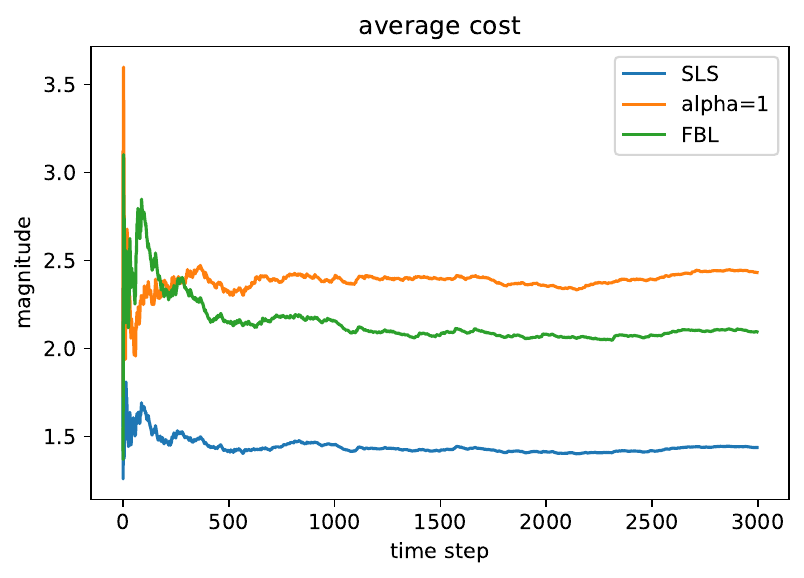}
    \caption{The time-averaged cost of the SLS controller is lower than the cost resulting from a traditional feedback linearization controller.}
    \label{fig:sls_cost}
\end{figure}

First, we fit a polynomial to the nonlinear function $f_a$ above, using order up to $k=3$ for the polynomial. We then generate the controller using a symbolic python package for an FIR horizon of $T=4$, which means that the effects of a disturbance arriving at time $t$ will be completely neutralized by time $t+4$. This results in 14 different functions of $\alpha_j^{(m)}$ to be computed at each time step. The NN has 8 linear layers of size 2048 each with ReLU and dropout activation functions between each layer. The stochastic gradient descent has a learning rate of 1e-6 and the disturbance bound is 1. The time-averaged cost of our controller is shown in Figure~\ref{fig:point_mass_cost}.
The feedback linearizing (FBL) controller is given by $u=g - \frac{1}{m} f_a(z_t,\dot{z}_t) + v_t$, where $v_t$ is selected to minimize the quadratic cost, as determined by the discrete Ricatti equation with $Q=I$ and $R=I$.

We see that our controller has a lower average cost compared with the FBL controller by about $33\%$. In the case where we set all $\alpha_j^{(k)}$ to one, which we test as a baseline for an untrained NN, we can see that the untrained SLS controller performs worse than the FBL controller. A key difference between the untrained SLS and FBL controllers is that the FBL controller makes the state exponentially decay, while the untrained SLS controller cancels disturbances in finite time, which are inherently different behaviors. We do not necessarily expect the untrained SLS controller to outperform the FBL controller because the FBL controller is optimal for the quadratic cost without disturbances.

\section{Conclusion and Outlook}\label{sec:conclusion}
We showed that using a polynomial approximation of nonlinear dynamics and applying the SLS controller leads to stable closed-loop dynamics, and that the control cost can be decreased by parameterizing the coefficients as NN functions. The stability guarantees hold regardless of the how well the NN is optimized, and when optimized well, we illustrate the tradeoffs among the stability size $h$, polynomial order $k$ and the disturbance size $W$. Future work will investigate the explicit relationship between the coefficients and the disturbances, formalize the controller paramterization for dynamics of the form $x_{t+1} = f(x_t) + g(x_t)u_t$, and implement this controller on a physical system.

\acks{Thanks to Guanya Shi, Eitan Levin, Chris Yeh, and Steven L. Brunton for helpful technical discussions around machine learning and polynomial systems. LC is supported by a NDSEG Graduate Fellowship from the Air Force Office of Scientific Research.}

\bibliography{papers}

\end{document}